\documentclass[12pt]{article}
\usepackage{amssymb,amsmath}
\usepackage{hyperref}
\usepackage{graphicx}
\usepackage{color}

\begin{document}

\title{\LARGE\bf A new asymptotic expansion series for the constant pi}

\author{
\normalsize\bf S. M. Abrarov\footnote{\scriptsize{Dept. Earth and Space Science and Engineering, York University, Toronto, Canada, M3J 1P3.}}\, and B. M. Quine$^{*}$\footnote{\scriptsize{Dept. Physics and Astronomy, York University, Toronto, Canada, M3J 1P3.}}}

\date{February 28, 2016}
\maketitle

\begin{abstract}
In our recent publications we have introduced the incomplete cosine expansion of the sinc function for efficient application in sampling [\href{http://dx.doi.org/10.1016/j.amc.2015.01.072}{Abrarov \& Quine, Appl. Math. Comput., 258 (2015) 425-435}; \href{http://dx.doi.org/10.5539/jmr.v7n2p163}{Abrarov \& Quine, J. Math. Research, 7 (2) (2015) 163-174}]. Here we show that it can also be utilized as a flexible and efficient tool in mathematical analysis. In particular, an application of the incomplete cosine expansion of the sinc function leads to expansion series of the error function in form of a sum of the Gaussian functions. This approach in integration provides a new asymptotic formula for the constant $\pi$.
\vspace{0.25cm}
\\
\noindent {\bf Keywords:} sinc function, error function, Gaussian function, asymptotic expansion, constant pi
\vspace{0.25cm}
\end{abstract}

\section{Introduction}
The sinc function is defined as \cite{Gearhart1990, Kac1959}
\begin{equation}\label{eq_1}
{\rm{sinc}}\left( x \right) = \left\{ 
\begin{aligned}
\frac{{\sin x}}{x}, \quad & x \ne 0\\
1, \quad & x = 0.
\end{aligned} \right.
\end{equation}
Due to its unique properties the sinc function is widely used in many applications in Applied Mathematics and Computational Physics. In particular, it can be applied in the Fourier analysis, numerical integration and differentiation \cite{Stenger2011}.

Vieta discovered a remarkable identity showing that the sinc function can be expressed as an infinite product \cite{Gearhart1990, Kac1959}
\begin{equation}\label{eq_2}
{\rm{sinc}}\left( x \right) = \prod\limits_{m = 1}^\infty  {\cos \left( {\frac{x}{{{2^m}}}} \right)}.
\end{equation}
The infinite product \eqref{eq_2} can be represented in form
\[
\prod\limits_{m = 1}^\infty  {\cos \left( {\frac{x}{{{2^m}}}} \right)}  = \left( {\prod\limits_{m = 1}^M {\cos \left( {\frac{x}{{{2^m}}}} \right)} } \right)\left( {\prod\limits_{n = 1}^\infty  {\cos \left( {\frac{x}{{{2^{M + n}}}}} \right)} } \right)
\]
and because of the limit
\[
\mathop {\lim }\limits_{M \to \infty } \prod\limits_{n = 1}^\infty  {\cos \left( {\frac{x}{{{2^{M + n}}}}} \right)}  = 1,
\]
we can infer that
\[
\prod\limits_{n = 1}^\infty  {\cos \left( {\frac{x}{{{2^{M + n}}}}} \right)}  \approx 1,		\qquad M > > 1.
\]
Consequently, the infinite product \eqref{eq_2} can be truncated as
$$
{\rm{sinc}}\left( x \right) \approx \prod\limits_{m = 1}^M {\cos \left( {\frac{x}{{{2^m}}}} \right)}.
$$
Although such a finite product of the cosine functions looks attractive for numerical integration, this form of the sinc approximation is not always convenient in practical applications since in many circumstances it is quite difficult to handle an integrand with product of the cosines. However, due to a product-to-sum identity that we have introduced in our previous publication \cite{Quine2013}
\begin{equation}\label{eq_3}
\prod\limits_{m = 1}^M {\cos \left( {\frac{x}{{{2^m}}}} \right)}  = \frac{1}{{{2^{M - 1}}}}\sum\limits_{m = 1}^{{2^{M - 1}}} {\cos \left( {\frac{{2m - 1}}{{{2^M}}}x} \right)},
\end{equation}
such a problem in truncation may be very effectively resolved. Specifically, substituting the product-to-sum identity \eqref{eq_3} into the Vieta\text{'}s infinite product \eqref{eq_2} yields the following limit
\begin{equation}\label{eq_4}
{\rm{sinc}}\left( x \right) = \mathop {\lim }\limits_{M \to \infty } \frac{1}{{{2^{M - 1}}}}\sum\limits_{m = 1}^{{2^{M - 1}}} {\cos \left( {\frac{{2m - 1}}{{{2^M}}}x} \right)}.
\end{equation}

The truncation of this equation to an integer ${2^{M - 1}}$ leads to the approximation based on the incomplete cosine expansion of the sinc function \cite{Abrarov2015a, Abrarov2015b}
\begin{equation}\label{eq_5}
{\rm{sinc}}\left( x \right) \approx \frac{1}{{{2^{M - 1}}}}\sum\limits_{m = 1}^{{2^{M - 1}}} {\cos \left( {\frac{{2m - 1}}{{{2^M}}}x} \right)}, \qquad - {2^{M - 1}}\pi  \le x \le {2^{M - 1}}\pi.
\end{equation}
Since the right side of this approximation is a sum of the cosine functions, it is periodic and, therefore, valid in the interval $ - {2^{M - 1}}\pi  \le x \le {2^{M - 1}}\pi $ near the origin.

As we have shown already, the upper limit in the incomplete cosine expansion of the sinc function may not be restricted by value ${2^{M - 1}}$ only \cite{Abrarov2015b}. It can also be extended to an arbitrary integer $L$. Such a transition from ${2^{M - 1}}$ to an arbitrary value $L$ can be shown by using the relation
$$
\cos \left( {\frac{{2m - 1}}{{{2^M}}}x} \right) = \cos \left( {\frac{{m - 1/2}}{{{2^{M - 1}}}}x} \right)
$$
and then by making change of the variable $L = {2^{M - 1}}$ in the equation \eqref{eq_4} that results in 
\begin{equation}\label{eq_6}
{\rm{sinc}}\left( x \right) = \mathop {\lim }\limits_{L \to \infty } \frac{1}{L}\sum\limits_{\ell  = 1}^L {\cos \left( {\frac{{\ell  - 1/2}}{L}x} \right)}.
\end{equation}
The presence of the limit \eqref{eq_6} signifies that if the integer $L$ is large enough, then it can be taken as an arbitrary integer leading to the following approximation \cite{Abrarov2015b} (see also equation (A.1) in the preprint version of Ref. \cite{Abrarov2015a})
\begin{equation}\label{eq_7}
{\rm{sinc}}\left( x \right) \approx \frac{1}{L}\sum\limits_{\ell  = 1}^L {\cos \left( {\frac{{\ell  - 1/2}}{L}x} \right)}, \qquad - \pi L \le x \le \pi L.
\end{equation}
Similar to equation \eqref{eq_5} this form of the incomplete cosine expansion of the sinc function is valid within the range $ - \pi L \le x \le \pi L$ as the right side of the equation \eqref{eq_7} is a periodic function. Thus, we can consider the approximation \eqref{eq_7} as a more generalized form of the incomplete cosine expansion \eqref{eq_5} of the sinc function. Particularly, comparing equations \eqref{eq_5} and \eqref{eq_7} with each other one can see that the approximation \eqref{eq_5} is a just specific case that occurs at $L = {2^{M - 1}}$.

The application of the incomplete cosine expansion \eqref{eq_7} of the sinc function is flexible and, therefore, efficient in numerical integration. In our previous publications we have shown that the incomplete cosine expansion of the sinc function can be effective in sampling \cite{Abrarov2015a, Abrarov2015b}. However, the sampling is not the only possibility where it can be very useful. In this work we show that application of the incomplete cosine expansion of the sinc function can also be used as a flexible tool for mathematical analysis. Specifically, using the incomplete cosine expansion of the sinc function we can derive a new asymptotic expansion series for the constant $\pi$.

\section{Results and discussion}

\subsection{The error function}

The most common form of definition of the error function is given by \cite{Abramowitz1972, Weisstein2003}
\begin{equation}\label{eq_8}
{\rm{erf}}\left( z \right) = \frac{2}{{\sqrt \pi  }}\int\limits_0^z {{e^{ - {t^2}}}dt},
\end{equation}
where $z = x + iy$ is the complex argument. The error function belongs to a family of special functions that sometimes regarded as the Faddeeva functions. This family of special functions includes the complex error function also known as the Faddeeva function \cite{Abrarov2015a, Faddeyeva1961}
$$
w\left( z \right) = {e^{ - {z^2}}}\left( {1 + \frac{{2i}}{{\sqrt \pi  }}\int\limits_0^z {{e^{{t^2}}}dt} } \right) = {e^{ - {z^2}}}\left[ {1 - {\rm{erf}}\left( { - iz} \right)} \right],
$$
the Dawson\text{'}s integral \cite{Weisstein2003, McCabe1974, Rybicki1989}
$$
{\rm{daw}}\left( z \right) = {e^{ - {z^2}}}\int\limits_0^z {{e^{t^2}}dt}  =  - \frac{{i\sqrt \pi  }}{2}{e^{ - {z^2}}}{\rm{erf}}\left( {iz} \right),
$$
the Fresnel integral \cite{Abramowitz1972}	
$$
F\left( z \right) = \int\limits_0^z {{e^{i\left( {\pi /2} \right){t^2}}}dt = \frac{{1 + i}}{2}{\rm{erf}}\left( {\frac{{\left( {1 - i} \right)\sqrt \pi  }}{2}z} \right)}
$$
and the normal distribution function \cite{Weisstein2003}
$$
\Phi \left( z \right) = \frac{1}{{\sqrt {2\pi } }}\int\limits_0^z {{e^{ - {t^2}/2}}dt = \frac{1}{2}{\rm{erf}}\left( {\frac{z}{{\sqrt 2 }}} \right)}.
$$

In general, the argument $z = x + iy$ of the error function is complex. However, in this paper we imply that its imaginary part $y$ of the argument $z$ is zero. Therefore, we only use the notation ${\rm{erf}}\left( x \right)$ signifying that $x \in \mathbb{R}$ is always real.

There is a remarkable integral of the error function (see integral 12 on page 4 in \cite{Ng1969})
\[
{\rm{erf}}\left( x \right) = \frac{1}{\pi }\int\limits_0^\infty  {{e^{ - u}}\sin \left( {2x\sqrt u } \right)\frac{{du}}{u}}.
\]
that can be readily expressed through the sinc function. In order to rearrange it accordingly, we make a change of the variable $v = \sqrt u$ leading to
\[
\begin{aligned}
{\rm{erf}}\left( x \right) &= \frac{1}{\pi }\int\limits_0^\infty  {{e^{ - {v^2}}}\sin \left( {2xv} \right)\frac{{2vdv}}{{{v^2}}}}  = \frac{2}{\pi }\int\limits_0^\infty  {{e^{ - {v^2}}}\sin \left( {2xv} \right)\frac{{dv}}{v}} \\
 &= \frac{{4x}}{\pi }\int\limits_0^\infty  {{e^{ - {v^2}}}\sin \left( {2xv} \right)\frac{{dv}}{{2xv}}} 
\end{aligned}
\]
or
\[
{\rm{erf}}\left( x \right) = \frac{{4x}}{\pi }\int\limits_0^\infty  {{e^{ - {v^2}}}{\rm{sinc}}\left( {2xv} \right)dv}.
\]
The factor $2$ that is present in the argument of the sinc function can be easily excluded by making another change of the variable $t = 2v$ providing
\[
{\rm{erf}}\left( x \right) = \frac{{4x}}{\pi }\int\limits_0^\infty  {{e^{ - {t^2}/4}}{\rm{sinc}}\left( {xt} \right)\frac{{dt}}{2}}
\]
or
\begin{equation}\label{eq_9}
{\rm{erf}}\left( x \right) = \frac{{2x}}{\pi }\int\limits_0^\infty  {{e^{ - {t^2}/4}}{\rm{sinc}}\left( {xt} \right)dt}.
\end{equation}
Further we will use this integral to approximate the error function in form of the Gaussian expansion series.

\subsection{Derivation}

There are many elegant expansions for the constant $\pi$ have been published in the modern literature \cite{Chudnovsky1989, Lange1999, Almkvist2003, Baruah2008, Baruah2009, Borwein2015}. Significant part of the recent developments in this field is related to the Ramanujan-type expansion series \cite{Baruah2008, Baruah2009, Borwein2015}. Nowadays, more than $12$ trillion digits of pi have been computed. Here we derive a new asymptotic expansion series for the constant $\pi$ by using the incomplete cosine expansion \eqref{eq_7} of the sinc function.

We may attempt to substitute the approximation \eqref{eq_7} for the sinc function into the integral \eqref{eq_9} in order to approximate the error function
\begin{equation}\label{eq_10}
{\rm{erf}}\left( x \right) \approx \frac{{2x}}{\pi }\int\limits_0^\infty  {\exp \left( { - {t^2}/4} \right)\underbrace {\frac{1}{L}\sum\limits_{\ell  = 1}^L {\cos \left( {\frac{{\ell  - 1/2}}{L}xt} \right)} }_{ \approx {\rm{sinc}}\left( {xt} \right)}dt}.
\end{equation}
However, since right side of the approximation \eqref{eq_7} is periodic, we have to justify that such a substitution is possible.

As it has been mentioned above, the sinc function \eqref{eq_1} can be approximated by the incomplete cosine expansion \eqref{eq_7} only within the range $ - \pi L \le x \le \pi L$ near the origin. This can be seen from the Fig. 1 showing the incomplete cosine expansion (black curve) and the original sinc function (shadowed blue curve). For the specific example $L = 15$ shown in the Fig. 1 the corresponding range where the incomplete cosine expansion coincides with the original sinc function is $ - {\rm{47}}{\rm{.1239}} \le x \le {\rm{47}}{\rm{.1239}}$ since $\pi  \times 15 \approx 47.1239$. As we can see from this figure, due to periodicity the incomplete cosine expansion \eqref{eq_7} cannot approximate the original sinc function beyond this specified range.

According to Milone {\it{et al.}} \cite{Milone1988}, an integral of kind
$$
\int\limits_{ - \infty }^\infty  {{e^{ - {t^2}}}f\left( t \right)dt},
$$
where $f\left( t \right)$ is a bounded function, can be approximated as
\begin{equation}\label{eq_11}
\int\limits_{ - \infty }^\infty  {{e^{ - {t^2}}}f\left( t \right)dt}  \approx \int\limits_{ - 6}^6 {{e^{ - {t^2}}}f\left( t \right)dt}
\end{equation}
retaining high accuracy since ${e^{ - {t^2}}}$ very rapidly decreases when $\left| t \right| > 6$ effectively damping such a way the entire integrand to zero. From the relation \eqref{eq_11} it immediately follows that for a bounded function $g\left( {xt} \right)$
$$
\int\limits_{ - \infty }^\infty  {{e^{ - {t^2}/4}}g\left( {xt} \right)dt}  \approx \int\limits_{ - 12}^{12} {{e^{ - {t^2}/4}}g\left( {xt} \right)dt}
$$
or
$$
\int\limits_0^\infty  {{e^{ - {t^2}/4}}g\left( {xt} \right)dt}  \approx \int\limits_0^{12} {{e^{ - {t^2}/4}}g\left( {xt} \right)dt}.
$$

\begin{figure}[ht]
\begin{center}
\includegraphics[width=22pc]{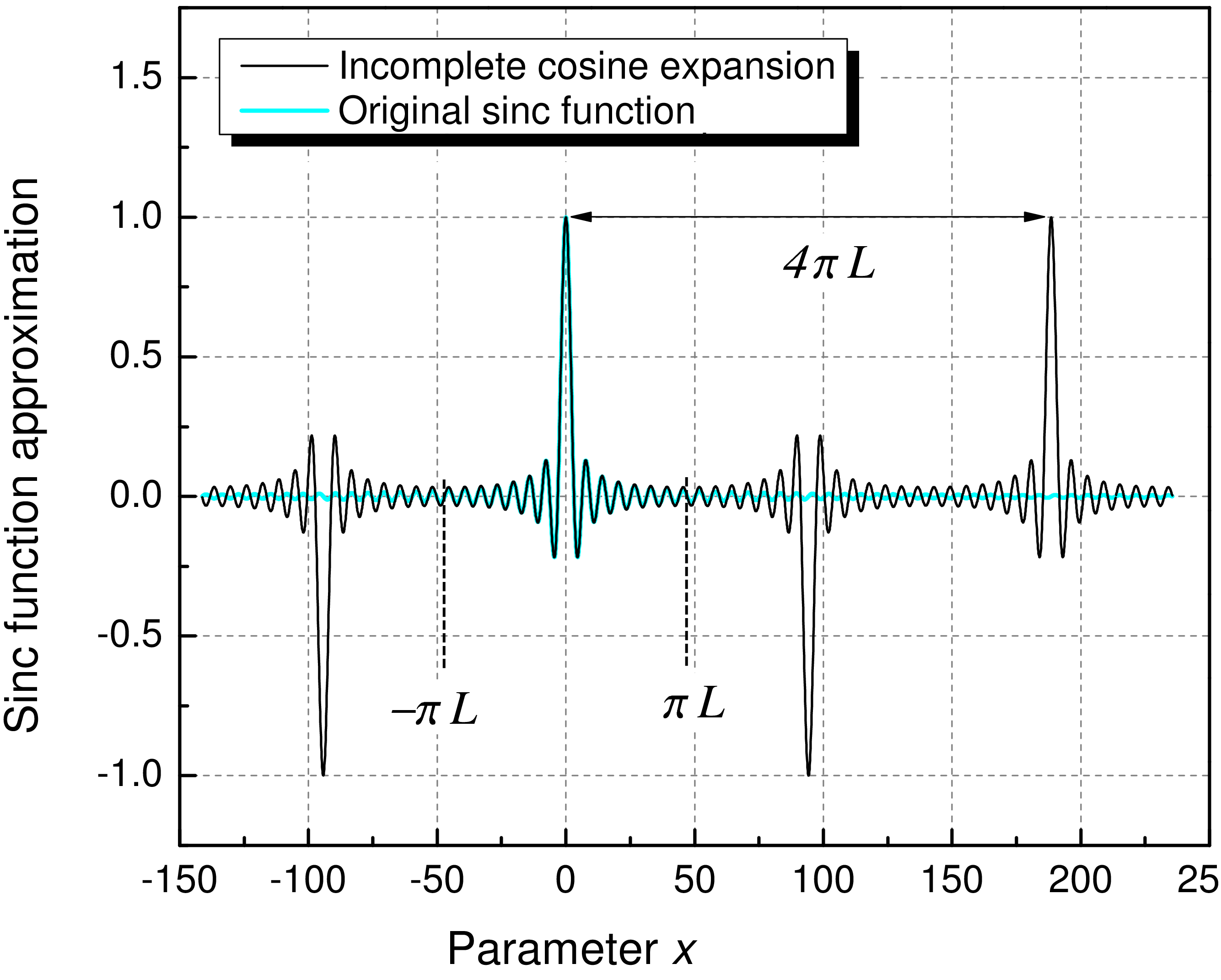}\hspace{2pc}%
\begin{minipage}[b]{28pc}
\vspace{0.5cm}
{\sffamily {\bf{Fig. 1.}} The incomplete cosine expansion (black curve) coincides with original the sinc function (shadowed blue curve) only inside the range $ - \pi L \le x \le \pi L$. At integer $L = 15$ the corresponding period is $4\pi L \approx {\rm{188}}{\rm{.4956}}$.}
\end{minipage}
\end{center}
\end{figure}

We can see now that despite periodicity, the approximation \eqref{eq_10} can be valid when the incomplete cosine expansion coincides with the original sinc function in the interval $0 \le x\,t \le 12$ since all other range is effectively damped to zero due to presence of the multiplier ${e^{ - {t^2}/4}}$ and, therefore, do not contribute for integration. Consequently, if the condition
\begin{equation}\label{eq_12}
\pi L \ge 12x
\end{equation}
is satisfied, then periodicity of the incomplete cosine expansion \eqref{eq_7} of the sinc function can be ignored in substitution into the approximation \eqref{eq_10}.

Each integral term on the right side of equation \eqref{eq_10} is integrable. Consequently, once the condition \eqref{eq_12} is satisfied, integration of the equation \eqref{eq_10} leads to the error function approximation in form of the Gaussian expansion series
\begin{equation}\label{eq_13}
{\rm{erf}}\left( x \right) \approx \frac{{2x}}{{\sqrt \pi  L}}\sum\limits_{\ell  = 1}^L {{e^{ - \frac{{{{\left( {\ell  - 1/2} \right)}^2}{x^2}}}{{{L^2}}}}}}  = \frac{{2x}}{{\sqrt \pi  L}}{\sum\limits_{\ell  = 1}^L {\left( {{e^{ - \frac{{{x^2}}}{{4{L^2}}}}}} \right)} ^{{{\left( {2\ell  - 1} \right)}^2}}}.
\end{equation}

Figure 2 shows some approximation curves for the error function obtained at $L = 5$, $L = 6$, $L = 7$, $L = 8$ by brown, green, red and blue curves, respectively. The black curve corresponding to the original error function is also shown for comparison. As we can see from this figure the curves tend to zero earlier and faster with decreasing $L$. This can be readily explained since the criterion \eqref{eq_12} is violated stronger at smaller values of the integer $L$ and larger values of the argument $x$.

\begin{figure}[ht]
\begin{center}
\includegraphics[width=22pc]{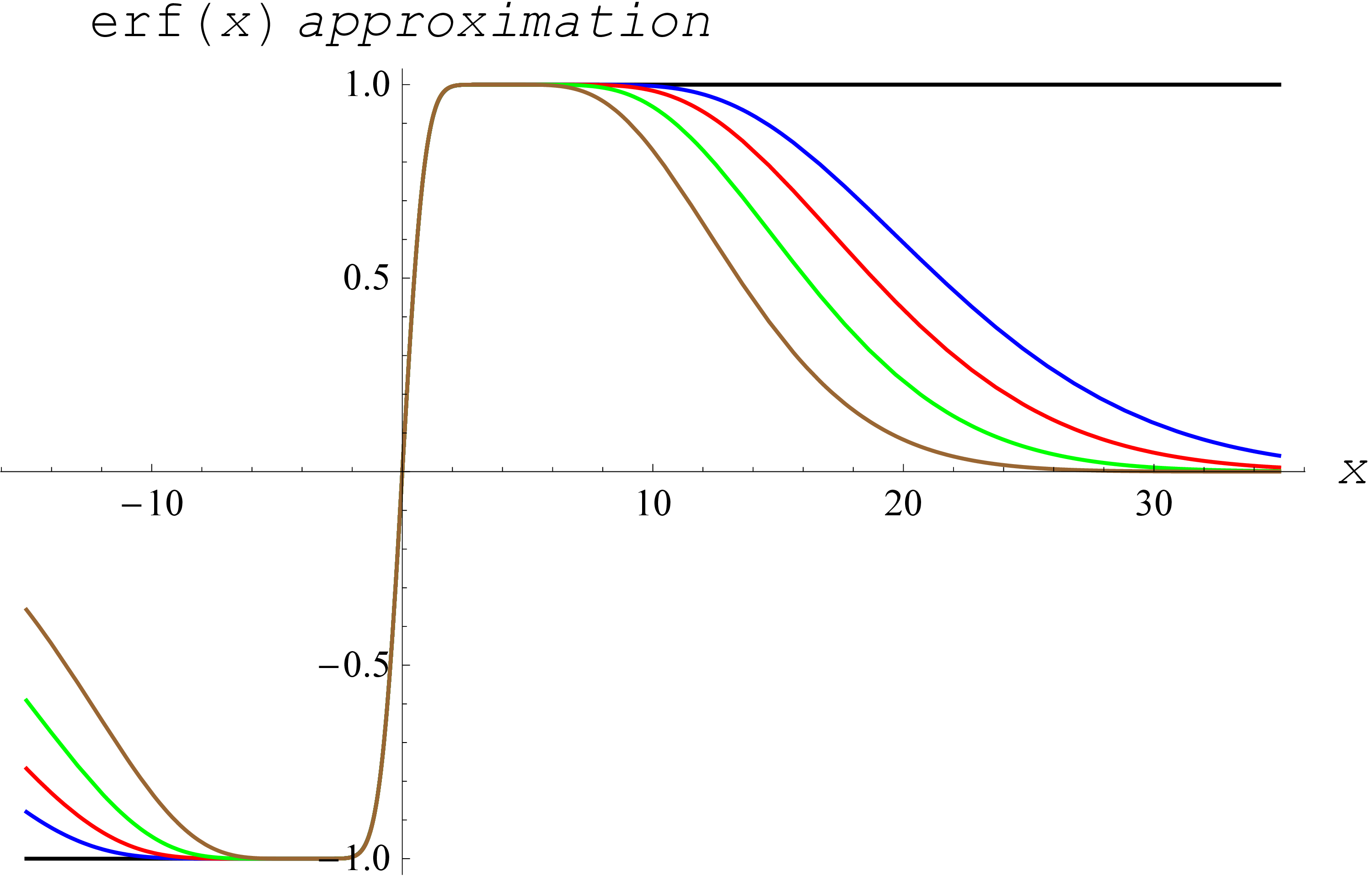}\hspace{2pc}%
\begin{minipage}[b]{28pc}
\vspace{0.75cm}
{\sffamily {\bf{Fig. 2.}} The error function approximations corresponding to $L = 5$ (brown curve), $L = 6$ (green curve), $L = 7$ (red curve) and $L = 8$ (blue curve). The original error function (black curve) is also shown for comparison.}
\end{minipage}
\end{center}
\end{figure}

It should be noted that because of the criterion \eqref{eq_12} the error function approximation \eqref{eq_13} can cover accurately only a central part near the origin. This problem, however, can be resolved since the computational test shows that ${\rm{erf}}\left( {x > 6} \right)$ and ${\rm{erf}}\left( {x < - 6} \right)$ are practically equal to $1$ and $- 1$, respectively. Consequently, we can cover the entire range $x \in \left( { - \infty ,\infty } \right)$ by using the following approximation
$$
{\rm{erf}}\left( x \right) \approx \left\{ 
\begin{aligned}
\frac{{2x}}{{\sqrt \pi  L}}\sum\limits_{\ell  = 1}^L {{e^{ - \frac{{{{\left( {\ell  - 1/2} \right)}^2}{x^2}}}{{{L^2}}}}}} , \qquad &  - 6 \le x \le 6\\
1, \qquad  & \,\, {\rm{otherwise}}.
\end{aligned} \right.
$$

Consider the following integral (see integral 1 on page 7 in \cite{Ng1969})
$$
\int\limits_0^a {{e^{ - {t^2}}}{\rm{erf}}\left( t \right)dt = \frac{{\sqrt \pi  }}{4}{{\left( {{\rm{erf}}\left( a \right)} \right)}^2}}.
$$
From the definition \eqref{eq_8} for the error function it is not difficult to see that ${\rm{erf}}\left( \infty  \right) = 1$. Consequently, by stretching $a$ to infinity we can rewrite this integral as an identity for the square root of pi
\begin{equation}\label{eq_14}
\sqrt \pi   = 4\int\limits_0^\infty  {{e^{ - {t^2}}}{\rm{erf}}\left( t \right)dt}.
\end{equation}
Since the error function is a bounded function $ - 1 \le {\rm{erf}}\left( x \right) \le 1$, 
then, according to Milone {\it{et al.}} \cite{Milone1988}, we can write
$$
\int\limits_0^\infty  {{e^{ - {t^2}}}{\rm{erf}}\left( t \right)dt}  \approx \int\limits_0^6 {{e^{ - {t^2}}}{\rm{erf}}\left( t \right)dt}
$$
and since the argument $t$ in the integrand varies from $0$ to $6$, the criterion \eqref{eq_12} for this case can be written as
$$ 
\pi L \ge 12 \times 6.
$$
Thus, we have reached an important result; we estimate a smallest value of the integer $L$ to be equal to $12 \times 6/\pi  \approx {\rm{23}}$ at which we can ignore periodicity of the incomplete cosine expansion \eqref{eq_7} of the sinc function. Assuming now that the integer $L \ge 23$, the substitution of the error function approximation \eqref{eq_13} into identity \eqref{eq_14} provides
$$
\sqrt \pi   \approx 4\int\limits_0^\infty  {\exp \left( { - {t^2}} \right)\underbrace {\frac{{2t}}{{\sqrt \pi  L}}\sum\limits_{\ell  = 1}^L {{e^{ - \frac{{{{\left( {\ell  - 1/2} \right)}^2}{t^2}}}{{{L^2}}}}}} }_{ \approx {\rm{erf}}\left( t \right)}dt}.
$$
In this equation each integral term is integrable. As a result, we have
$$
\sqrt \pi   \approx \frac{{16}}{{\sqrt \pi  }}L\sum\limits_{\ell  = 1}^L {\frac{1}{{{{\left( {2\ell  - 1} \right)}^2} + 4{L^2}}}} 
$$
or
\begin{equation}\label{eq_15}
\pi  \approx 16L\sum\limits_{\ell  = 1}^L {\frac{1}{{{{\left( {2\ell  - 1} \right)}^2} + 4{L^2}}}}  = 16L\sum\limits_{\ell  = 1}^L {\frac{1}{{4{\ell ^2} - 4\ell  + 1 + 4{L^2}}}}.
\end{equation}
Since ${L^2} >  > 1$, we can simply ignore $ - 4\ell  + 1$ in the denominator due to vanishing contribution. Consequently, we obtain the asymptotic expansion series for the constant $\pi$
\begin{equation}\label{eq_16}
\pi  \approx 4L\sum\limits_{\ell  = 1}^L {\frac{1}{{{\ell ^2} + {L^2}}}}.
\end{equation}

The expansion series \eqref{eq_15} and \eqref{eq_16} are asymptotic because they both converge to the constant $\pi $ as the integer $L$ approaches to infinity. To the best of our knowledge the asymptotic formulas \eqref{eq_15} and \eqref{eq_16} for computing the constant $\pi $ have never been reported in scientific literature. Perhaps these asymptotic formulas may be grouped as a special kind of expansions due to their specific feature - the presence of an arbitrarily large integer $L >> 1$ that is involved in each summation term in computing pi.

There is an easier way to derive the asymptotic series \eqref{eq_15} and \eqref{eq_16}. Since the equation \eqref{eq_6} is exact, its substitution into identity \eqref{eq_9} yields
\begin{equation}\label{eq_17}
{\rm{erf}}\left( x \right) = \frac{{2x}}{\pi } \times \mathop {\lim }\limits_{L \to \infty } \int\limits_0^\infty  {\exp \left( { - {t^2}/4} \right)\underbrace {\frac{1}{L}\sum\limits_{\ell  = 1}^L {\cos \left( {\frac{{\ell  - 1/2}}{L}xt} \right)} }_{{\rm{sinc}}\left( {xt} \right)}dt}
\end{equation}
or
\begin{equation}\label{eq_18}
{\rm{erf}}\left( x \right) = \frac{{2x}}{{\sqrt \pi  }} \times \mathop {\lim }\limits_{L \to \infty } \frac{1}{L}\sum\limits_{\ell  = 1}^L {{e^{ - \frac{{{{\left( {\ell  - 1/2} \right)}^2}{x^2}}}{{{L^2}}}}}}.
\end{equation}
Once again, since the equation \eqref{eq_18} is exact its substitution into identity \eqref{eq_14} provides
$$
\sqrt \pi   = 4 \times \mathop {\lim }\limits_{L \to \infty } \int\limits_0^\infty  {\exp \left( { - {t^2}} \right)\underbrace {\frac{{2t}}{{\sqrt \pi  L}}\sum\limits_{\ell  = 1}^L {{e^{ - \frac{{{{\left( {\ell  - 1/2} \right)}^2}{t^2}}}{{{L^2}}}}}} }_{{\rm{erf}}\left( t \right)}dt}
$$
or
$$
\sqrt \pi   = \frac{{16}}{{\sqrt \pi  }} \times \mathop {\lim }\limits_{L \to \infty } \sum\limits_{\ell  = 1}^L {\frac{L}{{{{\left( {2\ell  - 1} \right)}^2} + 4{L^2}}}}
$$
or
\begin{equation}\label{eq_19}
\pi  = 16 \times \mathop {\lim }\limits_{L \to \infty } \sum\limits_{\ell  = 1}^L {\frac{L}{{{{\left( {2\ell  - 1} \right)}^2} + 4{L^2}}}}  = 16 \times \mathop {\lim }\limits_{L \to \infty } \sum\limits_{\ell  = 1}^L {\frac{L}{{4{\ell ^2} - 4\ell  + 1 + 4{L^2}}}}.
\end{equation}
While $L$ tends to infinity ${L^2} >  > 1$ and ${L^2} >  > \ell $. Consequently, we can simplify the limit \eqref{eq_19} as given by
\begin{equation}\label{eq_20}
\pi  = 4 \times \mathop {\lim }\limits_{L \to \infty } \sum\limits_{\ell  = 1}^L {\frac{L}{{{\ell ^2} + {L^2}}}}.
\end{equation}

Truncating now the equations \eqref{eq_19} and \eqref{eq_20} we obtain the asymptotic expansion series \eqref{eq_15} and \eqref{eq_16}, respectively, for the constant pi. Although this way of derivation is straightforward, it does not explain why the limits \eqref{eq_19} and \eqref{eq_20} can be truncated for computation of $\pi $ if the finite value of the integer $L$ signifies an appearance of the periodicity on right side of the equation \eqref{eq_17}. Furthermore, this derivation also does not provide an estimated value for a smallest $L$ at which the periodicity can be ignored. Therefore, the initial derivation discussed earlier in this section is more informative and preferable.

\subsection{Sample computations}

In this section we represent some examples of approximated constant pi and compare them with the reference
$$
\pi  = 3.1415926535897932384626433832795028841 \ldots \,\,.
$$
We performed all computational tests by using Wolfram Mathematica 9 in enhanced precision mode in order to demonstrate how many digits can coincide with the reference at some given $L$.

As we have found the estimated value for the smallest integer, it would be reasonable to start computation of the constant pi at $L = 23$. Substitution shows that even with this small integer the equation \eqref{eq_15} provides a number where four digits coincide with the reference
$$
\pi  \approx \underbrace {3.141}_{{\rm{4}}\,\,{\rm{digits}}}7501835074959226730514058333 \ldots \,\,.
$$
However, with same integer the equation \eqref{eq_16} cannot compute properly. In particular, the result of computation is
$$
\pi  \approx \underbrace 3_{{\rm{1}}\,\,{\rm{digit}}}.0977993328722573504830557765576 \ldots \,\,,
$$
where only first digit coincides with the reference. Such a low accuracy occurs due to simplification we made by excluding $ - 4\ell  + 1$ terms in denominator of the equation \eqref{eq_16} under assumption that $L$ is large.  As a result, when $L$ is relatively small, say below ${10^4}$, we should not expect a close approximation for pi.

Increase of the integer up to $L = {10^3}$ shows the values
$$
\pi  \approx \underbrace {3.141592}_{{\rm{7}}\,\,{\rm{digits}}}736923126571794054593597 \ldots
$$
and
$$
\pi  \approx \underbrace {3.14}_{{\rm{3}}\,\,{\rm{digits}}}0\underbrace {592}_{{\rm{3}}\,\,{\rm{digits}}}486923126571797960843597 \ldots \,\,
$$
for the equations \eqref{eq_15} and \eqref{eq_16}, respectively. The equation \eqref{eq_16} shows unusual behavior; the digits coinciding with the reference consist of two groups of three digits separated by a single non-coinciding digit $0$. Furthermore, the non-coinciding digit $0$ is smaller than the actual digit $1$ by one.

Increase of the integer up to $L = {10^6}$ shows a significant improvement in accuracy of the equation \eqref{eq_15}
$$
\pi  \approx \underbrace {3.141592653589}_{13\,\,{\rm{digits}}}87657179597671661 \ldots \,\,.
$$
As we can see, the thirteen digits coincide with the reference. However, with same integer, the equation \eqref{eq_16} again shows unusual behavior; the digits coinciding with the reference consist of two groups of six digits separated by a single non-coinciding digit.
$$
\pi  \approx \underbrace {3.14159}_{{\rm{6}}\,\,{\rm{digits}}}1\underbrace {653589}_{{\rm{6}}\,\,{\rm{digits}}}62657179597671661 \ldots \,\,.
$$
We also note that this non-coinciding digit $1$ between two groups is smaller than the corresponding actual digit $2$ by one. 

Further increase of the integer up to $L = {10^9}$ in the equation \eqref{eq_15} yields further improvement in accuracy revealing $19$ coinciding digits with the reference
$$
\pi  \approx \underbrace {3.141592653589793238}_{{\rm{19}}\,\,{\rm{digits}}}54597671661 \ldots \,\,.
$$
However, the equation \eqref{eq_16} again retain the unusual behavior in computation
$$
\pi  \approx \underbrace {3.14159265}_{{\rm{9}}\,\,{\rm{digits}}}2\underbrace {589793238}_{{\rm{9}}\,\,{\rm{digits}}}29597671661 \ldots \,\,.
$$
Specifically, it yields a value consisting of two groups of nine coinciding digits where a single non-coinciding digit $2$ separating these two groups is smaller than the actual digit by one.

Lastly, the computational test we performed at $L = {10^{12}}$  reveals a rapid convergence in the equation \eqref{eq_15}
$$
\pi  \approx \underbrace {3.141592653589793238462643}_{{\rm{25}}\,\,{\rm{digits}}}46661 \ldots
$$
since there are $25$ digits that coincide with the reference. Once again, the equation \eqref{eq_16} persistently retain unusual tendency in computation
$$
\pi  \approx \underbrace {3.14159265358}_{{\rm{12}}\,\,{\rm{digits}}}8\underbrace {793238462643}_{{\rm{12}}\,\,{\rm{digits}}}21661 \ldots
$$
showing two groups of coinciding digits with $12$ digits in each group. These two groups of coinciding digits are separated by non-coinciding digit $8$ that is smaller than the actual digit $9$ by one.

Although the equation \eqref{eq_15} is not as simple as the equation \eqref{eq_16}, the computational test we performed shows that it is essentially more rapid in convergence. Consequently, the asymptotic expansion series \eqref{eq_15} is the main result of this work. This application demonstrates that the incomplete cosine expansion of the sinc function that we have introduced earlier for sampling \cite{Abrarov2015a, Abrarov2015b}, can also be used as a flexible and efficient tool in mathematical analysis.

\section{Conclusion}
We show that the incomplete cosine expansion of the sinc function can be used as a flexible and efficient tool in mathematical analysis. Specifically, its application for the error function leads to the expansion series consisting of a sum of the Gaussian functions. Such an approach in integration provides a new asymptotic formula for the constant $\pi$. The computational test we performed shows a rapid convergence of the proposed asymptotic formula \eqref{eq_15}. In particular, at integer $L$ equal to $10^3$, $10^6$, $10^9$ and $10^{12}$ it provides respectively $7$, $13$, $19$ and $25$ coinciding digits with the reference value of $\pi$.

\section*{Acknowledgments}

This work is supported by National Research Council Canada, Thoth Technology Inc. and York University.


\end{document}